\documentclass[12pt]{article}

\usepackage{bm}
\usepackage{graphicx}
\usepackage{subfigure}
\usepackage{diagbox}
\usepackage{amsmath}
\usepackage{amssymb}
\usepackage{amsfonts}
\usepackage{fdsymbol}
\usepackage{hyperref}
\usepackage{color}
\hypersetup{colorlinks=true,linkcolor=blue,filecolor=blue,urlcolor=blue}

\title{Visualization for Dichotomous Variables,\\the Independence and Markov chains}
\author{ZHANG Yan\\ianzhang@connect.hku.hk}

\begin{document}

\maketitle

\section{Introduction}

In probability theory, the independence is a very fundamental concept, but with a little mystery. People can always easily manipulate it logistically but not geometrically, especially when it comes to the independence relationships among more that two variables, which may also involve conditional independence. Here I am particularly interested in visualizing Markov chains which have the well known memoryless property. I am not talking about drawing the transition graph, instead, I will draw all events of the Markov process in a single plot. Here, to simplify the question, this work will only consider dichotomous variables, but all the methods actually can be generalized to arbitrary set of discrete variables. Let me end the introduction with one final work below: 

\begin{figure}[ht!]
\centering
\includegraphics[width=0.5\linewidth]{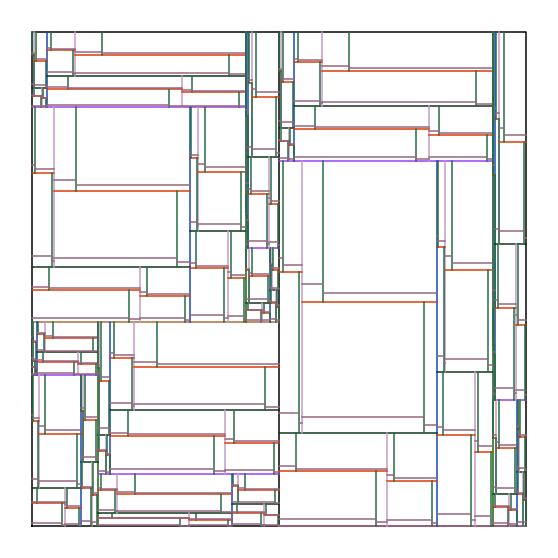}
\end{figure}
This is a non-time-homogeneous Markov chain with 10 dichotomous variables, and the conditional probabilities are randomly generated. We can actually observe some sort of fractal pattern in this graph, which is just the reflection of the Markov property. 

\section{Visualization for Dichotomous Variables}

Assume probability space $(\Omega,\mathcal{F},P)$, a dichotomous variable is just a function $X$: $\Omega\mapsto\{0,1\}$, which means a dichotomous variable can split the sample space $\Omega$ into two parts, $X^{-1}(0)$ and $X^{-1}(1)$. So, if I assume the points in a plane geometry like a circle, a triangle or a rectangle as $\Omega$, then I can use anther plane geometry or simply a line to represent a dichotomous variable if it split the $\Omega$ into two parts. Like what we can see below:

\begin{figure}[ht!]
\centering
\includegraphics[width=0.8\linewidth]{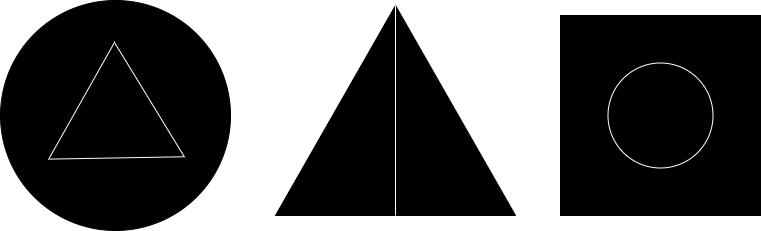}
\end{figure}
And actually, we don't set the constraint that the two parts should be two intact subspaces, which means the follow dichotomous variables are also valid, but in these cases, we need to set a rule to specify which subspaces are parts of $X^{-1}(0)$ and which are parts of $X^{-1}(1)$. 

\begin{figure}[ht!]
\centering
\includegraphics[width=0.8\linewidth]{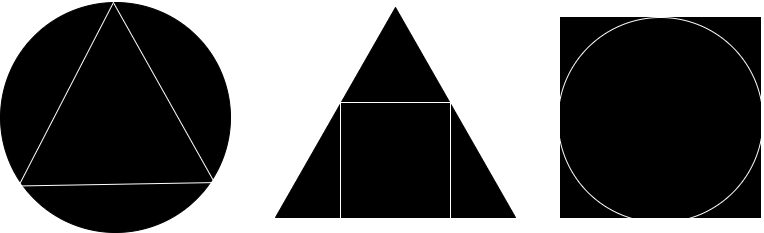}
\end{figure}
For the above example, we can assume the outer ones are parts of $X^{-1}(0)$, the inner one is $X^{-1}(1)$. I will avoid this kind of phenomenon if possible, but in the visualization of long Markov chain, it's not avoidable. 

\section{Visualization for the Independence}

With the assistance of squares, circles and lines, we can already visualize all kinds of ternary independence-dependence relationships quite elegantly.

\newpage
Firstly, I visualize the common effect relationship: 

\begin{figure}[ht!]
\centering
\includegraphics[width=0.6\linewidth]{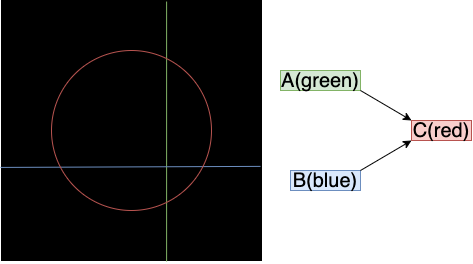}
\end{figure}
In the above graph, we can see that $A\nVbar C$, $B\nVbar C$, $A\Vbar B$ but $A\nVbar B|C$, which is the selection bias phenomenon in causal inference. 

And then, the common cause relationship:

\begin{figure}[ht!]
\centering
\includegraphics[width=0.6\linewidth]{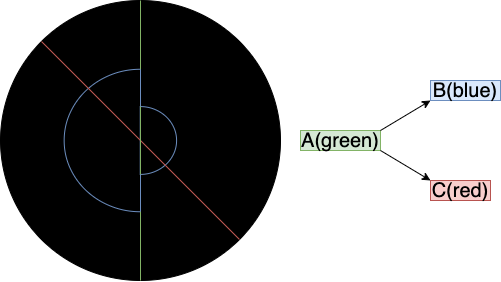}
\end{figure}
This time, $A\nVbar B$, $A\nVbar C$, $B\Vbar C|A$ but $B\nVbar C$, which is the confounding phenomenon in causal inference. And furthermore, this graph can also represent a Markov chain of length three: $B\to A\to C$, because $B\nVbar A$, $A\nVbar C$, $B\nVbar C$ but $B\Vbar C|A$. So, for this three variables, we actually can't tell the difference between common cause relationship and causal chain relationship. And by symmetry, the above graph can also be $C\to A\to  B$, which is a well-known phenomenon (Markov equivalence graphs \cite{meek2013causal}) that the following causal diagrams are equivalent to each other in a sense. 

\begin{figure}[ht!]
\centering
\includegraphics[width=0.6\linewidth]{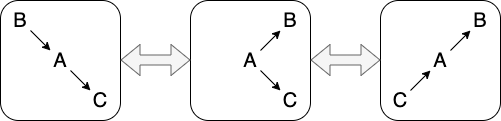}
\end{figure}

\section{Visualization for Markov chains}

To draw the Markov chains, we need to adopt some more complicated rules. Firstly, for simplicity, I set our aim to be generating a Markov chain $X_0\to X_1\to\dots\to X_n$ with symmetric time-homogeneous conditional probabilities: $P(X_{i+1}=0|X_{i}=0)=1/3$, $P(X_{i+1}=0|X_{i}=1)=2/3$, and the initial probability $P(X_0=1)=1/2$. Let me show you the first three graphs:

\begin{figure}[ht!]
\centering
\subfigure[$X_0$]{
\includegraphics[width=0.3\linewidth]{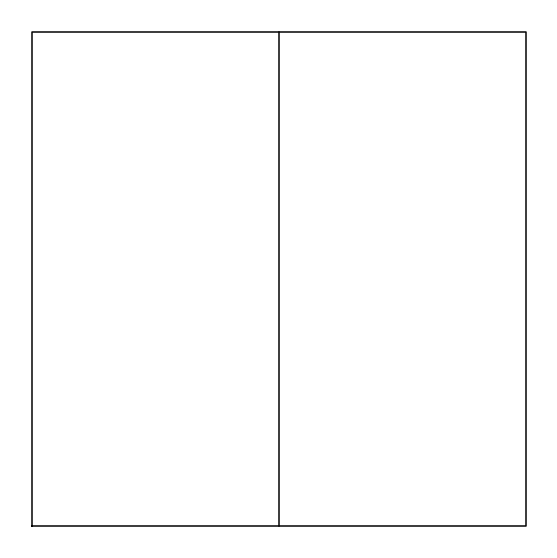}}
\subfigure[$X_0\to X_1$]{
\includegraphics[width=0.3\linewidth]{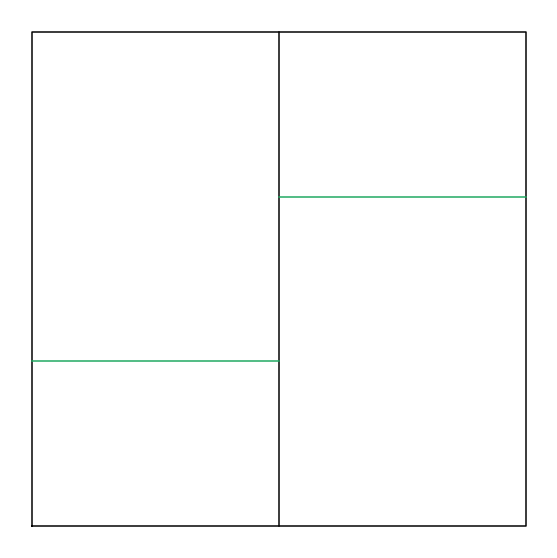}}
\subfigure[$X_0\to X_1\to X_2$]{
\includegraphics[width=0.3\linewidth]{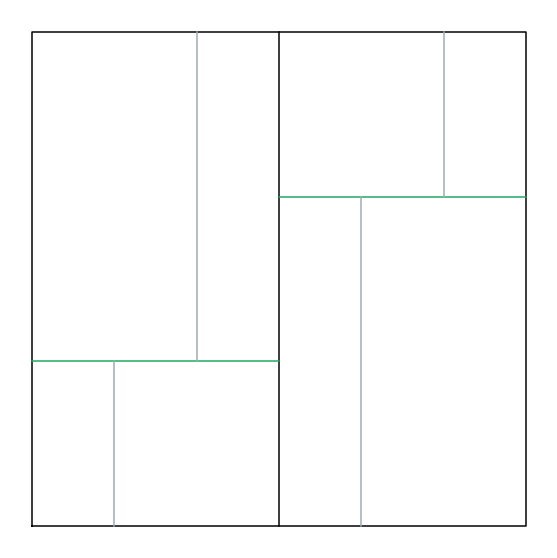}}
\end{figure}

In (a), one single vertical line was used to represent $X_0$, splitting the whole space into two equal parts, where the left one is $X_0^{-1}(0)$ and the right one is $X_0^{-1}(1)$. 

In (b), based on $P(X_1=0|X_0=0)=1/3$, I used a horizontal line to split subspace $X_0^{-1}(0)$ into two parts, the lower one $X_1^{-1}(0)\cap X_0^{-1}(0)$ and the upper one $X_1^{-1}(1)\cap X_0^{-1}(0)$. Similarly, by $P(X_1=0|X_0=1)=2/3$, two more parts, $X_1^{-1}(0)\cap X_0^{-1}(1)$ and $X_1^{-1}(1)\cap X_0^{-1}(1)$ were obtained. And $X_1^{-1}(0)$ is just the union of $X_1^{-1}(0)\cap X_0^{-1}(0)$ and $X_1^{-1}(0)\cap X_0^{-1}(1)$, and $X_1^{-1}(1)=(X_1^{-1}(1)\cap X_0^{-1}(0))\cup(X_1^{-1}(1)\cap X_0^{-1}(1))$. We can see in this plot, $X_0\nVbar X_1$. 

In (c), by $P(X_2=0|X_1=0)=1/3$ and $P(X_2=0|X_1=1)=2/3$, there are four subspaces for $X_2^{-1}(0)$ and $X_2^{-1}(1)$ formed respectively by four new vertical lines. Easy to see that $X_1\nVbar X_2$, and because $P(X_2=0|X_0=0)=5/9$ which is not equal to $P(X_2=0|X_0=1)=4/9$, $X_0\nVbar X_2$ can also be deduced. But this time, we actually have $X_0\Vbar X_2|X_1$, which makes the graph a Markov chain. 

In next page, I drew the plots from 2-chain to 10-chain, each color representing one single dichotomous variable. 

\begin{figure}[ht!]
\centering
\includegraphics[width=1\linewidth]{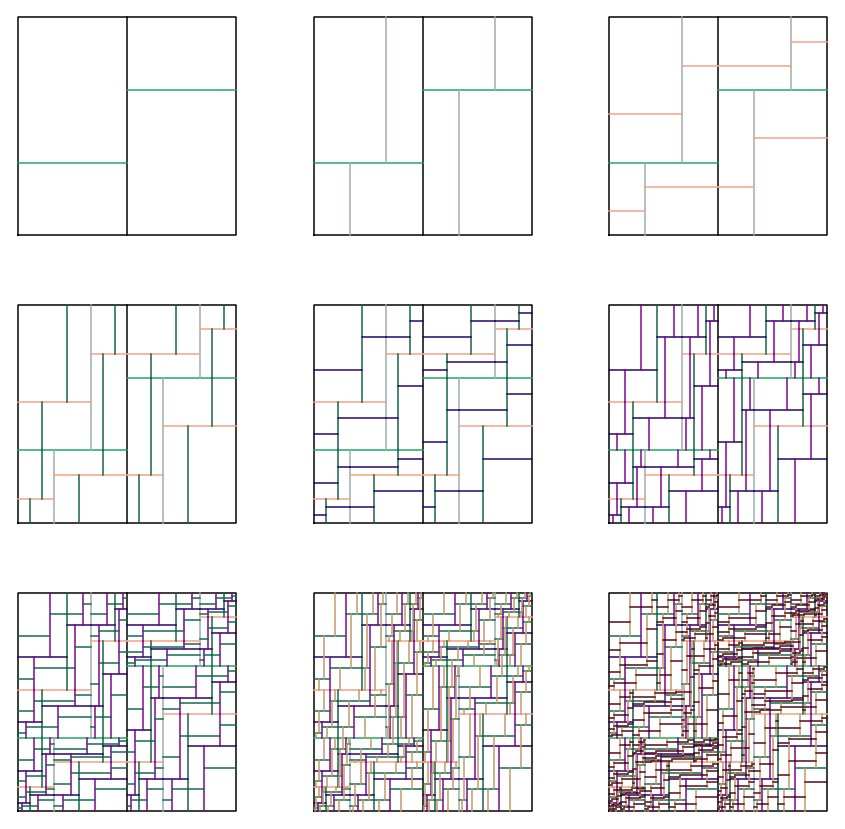}
\end{figure}

Finally, if I change the conditional probabilities and the initial probability, I can actually draw all kinds of Markov chains made by dichotomous variables as I like, such as the graph I showed in the introduction. And my current method is not the only way to visualize such complicated independence relationships, one can easily find another one if he understand the fundamental, and sometime it's just a matter of aesthetics.

\bibliographystyle{plain}
\bibliography{bib}

\end{document}